\documentclass[11pt,epsfig,amssymb,amsmath]{amsart}
\usepackage{amsmath,amssymb}
\makeatletter \oddsidemargin.9375in \evensidemargin \oddsidemargin
\begin{document}
\begin{center}
 {\large\bf Approximately  Lie ternary
$(\sigma,\tau,\xi)-$derivations on Banach ternary algebras
\vskip.20in

{\bf M. Eshaghi Gordji } \\[2mm]

{\footnotesize Department of Mathematics,
Semnan University,\\ P. O. Box 35195-363, Semnan, Iran\\
[-1mm] e-mail: {\tt madjid.eshaghi@gmail.com}}

{\bf R. Farrokhzad} \\[2mm]

{\footnotesize Department of Mathematics,
Shahid Beheshti University, Tehran, Iran\\
[-1mm] e-mail: {\tt razieh.farokhzad@yahoo.com}}

{\bf S. A. R. Hosseinioun} \\[2mm]

{\footnotesize Department of Mathematics,
Shahid Beheshti University, Tehran, Iran\\
[-1mm] e-mail: {\tt ahosseinioun@yahoo.com}}}

\vskip 5mm {\bf Abstract}
\end{center}
\vskip 5mm \noindent{\footnotesize \  \ \ Let $A$ be a Banach
ternary algebra over a scalar field $\Bbb R$ or $\Bbb C$ and $X$ be
a ternary Banach $A-$module. Let $\sigma,\tau$ and $\xi$ be linear
mappings on $A$, a linear mapping $D:(A,[~]_A)\to (X,[~]_X)$ is
called a Lie ternary $(\sigma,\tau,\xi)-$derivation, if
$$D([abc]_A)=[[D(a)bc]_X]_{(\sigma,\tau,\xi)}+[[D(b)ac]_X]_{(\sigma,\tau,\xi)}+[[D(c)ba]_X]_{(\sigma,\tau,\xi)},$$
for all
$a,b,c\in A$, where $[abc]_{(\sigma,\tau,\xi)}=a\tau(b)\xi(c)-\sigma(c)\tau(b)a.$\\
In this paper, we investigate  the generalized Hyers--Ulam--Rassias
stability of Lie ternary $(\sigma,\tau,\xi)-$derivations on Banach
ternary algebras.

\vskip.10in
 \footnotetext { 2000 Mathematics Subject Classification: 39B82;
 39B52; 46B99; 17A40.}
 \footnotetext { Keywords: Hyers--Ulam--Rssias stability; ternary algebra; Lie derivation}

\vskip .10in
\newtheorem{df}{Definition}[section]
\newtheorem{rk}[df]{Remark}
\newtheorem{lem}[df]{Lemma}
\newtheorem{thm}[df]{Theorem}
\newtheorem{pro}[df]{Proposition}
\newtheorem{cor}[df]{Corollary}
\newtheorem{ex}[df]{Example}
\setcounter{section}{0} \numberwithin{equation}{section} \vskip
3mm
\begin{center}
{\section{\bf Introduction}}
\end{center}\vskip 2mm
In the 19 th century, many mathematicians considered ternary algebraic operations and their generalizations.
 A.Cayley ([7])  introduced the notion of cubic matrix. It was later generalized by Kapranov,
 Gelfand and Zelevinskii in 1990 ([12]).
Below, a composition rule includes a simple example of such non-trivial ternary operation:
$$\{a,b,c\}_{ijk}=\sum_{l,m,n}a_{nil}b_{ljm}c_{mkn},\ \ \ \ \ \ \ \ \ \ {i,j,k...=1,2,...,N.}$$
There are a lot of hopes that ternary structures and their generalization will have certain possible applications
in physics. some of these applications are (see [2,3],[5],[10],[13,14,15]).
A ternary (associative) algebra ($A$, [ ]) is a linear space $A$ over a scalar
field $\Bbb F=(\Bbb R$ or $\Bbb C)$ equipped with a linear mapping,
 the so-called ternary product, [ ]: $A\times A\times A\to A$ such that
 $[[abc]de]=[a[bcd]e]$ for all $a,b,c,d,e \in A$. This notion is a natural
 generalization of the binary case. Indeed if $(A,\odot)$ is a usual (binary) algebra
 then $[abc]:=(a\odot b)\odot c$ induced a ternary product making $A$ into a
ternary algebra which will be called trivial. It is known that unital ternary algebras
are trivial and finitely generated ternary algebras are ternary subalgebras of trivial ternary algebras [6].
There are other types of ternary algebras in which one may consider other
versions of associativity. Some examples of ternary algebras are $(\textit{i})$ "cubic matrices" introduced
by Cayley [7] which were in turn generalized by Kapranov, Gelfand and Zelevinskii [12];
$(\textit{ii})$ the ternary algebra of polynomials of odd degrees in one variable
equipped with the ternary operation $[p_1p_2p_3]=p_1\odot p_2\odot p_3$,
where $\odot$ denotes the usual multiplication of polynomials.\\
\\ By a Banach ternary algebra we mean a ternary algebra equipped with a complete norm $\|.\|$ such that
 $\|[abc]\|\leq \|a\|\|b\|\|c\|$.
If a ternary algebra ($A$,[ ]) has an identity, i.e. an element $e$ such that $a=[aee]=[eae]=[eea]$ for all
$a\in A$, then $a\odot b:=[aeb]$ is a binary product for which we have
$$(a\odot b)\odot c=[[aeb]ec]=[ae[bec]]=a\odot(b\odot c)$$
and
$$a\odot e=[aee]=a=[eea]=e\odot a,$$
for all $a,b,c\in A$ and so ($A$,[ ]) may be considered as a (binary) algebra.
Conversely, if $(A,\odot)$ is any (binary) algebra, then $[abc]:=a\odot b\odot c$ makes $A$
into a ternary algebra with the unit $e$ such that $a\odot b=[aeb]$.\\
\\ Let $A$ be a Banach ternary algebra and $X$ be a Banach space.
Then $X$ is called a ternary Banach $A-$module, if module operations
$A \times A \times X \to X,$ $A \times X \times A \to X,$ and $X
\times A \times A \to X$ are $\Bbb C-$linear in every variable.
Moreover  satisfy:
$$[[abc]_A~dx]_X=[a[bcd]_A~x]_X=[ab[cdx]_X]_X$$
$$[abc]_A~xd]_X=[a[bcx]_X~d]_X=[ab[cxd]_X]_X,$$
$$[[xab]_X~cd]_X=[x[abc]_A~d]_X=[xa[bcd]_A]_X,$$
$$[[axb]_X~cd]_X=[a[xbc]_X~d]_X=[ax[bcd]_A]_X,$$
$$[[abx]_X~cd]_X=[a[bxc]_X~d]_X=[ab[xcd]_X]_X,$$
for all $x \in X$ and all $a,b,c,d \in A,$ and
$$\max\{\|[xab]_X\|,\|[axb]_X\|,\|[abx]_X\|\}\leq \|a\| \|b\| \|x\|$$
for all $x \in X$ and all $a,b \in A.$\\
Let $A$ be a normed algebra, $\sigma$ and $\tau$ two mappings on $A$
and $X$ be an $A-$bimodule. A linear mapping $L:A\to X$ is called a
Lie $(\sigma,\tau)-$derivation, if
$$L([a,b])=[L(a),b]_{\sigma,\tau}-[L(b),a]_{\sigma,\tau}$$
for all $a,b\in A$, where $[a,b]_{\sigma,\tau}$ is $a\tau(b)-\sigma(b)a$ and $[a,b]$
is the commutator $ab-ba$ of elements $a,b.$\\
Now, Let $(A,[~]_A)$ be a Banach ternary algebra over a scalar field
$\Bbb R$ or $\Bbb C$ and $(X,[~]_X)$ be a ternary Banach $A-$module.
Let $\sigma, \tau$ and $\xi$ be linear mappings on $A$. A linear
mapping $D:(A,[~]_A) \to (X,[~]_X)$ is called a Lie ternary
$(\sigma,\tau,\xi)-$derivation, if
$$D([abc]_A)=[[D(a)bc]_X]_{(\sigma,\tau,\xi)}+[[D(b)ac]_X]_{(\sigma,\tau,\xi)}+[[D(c)ba]_X]_{(\sigma,\tau,\xi)}\eqno(1.1)$$ for all
$a,b,c\in A$, where $[abc]_{(\sigma,\tau,\xi)}=a\tau(b)\xi(c)-\sigma(c)\tau(b)a.$\\
\\
If a Banach ternary algebra $A$ has an identity $e$ such that
$\|e\|=1$, as we said above, $A$ may be considered as a (binary)
algebra. Now let $X$ be a ternary Banach $A-$module, then $X$ may be
considered as a Banach $A-$module by following module product:
$$a.x=[aex]_X\hspace{2.5cm}x.a=[xea]_X$$
for all $a\in A, x\in X.$\\
Let $A$ be a unital Banach ternary algebra and $X$ be a ternary
Banach $A-$module. If $D:A\to X$ is a Lie ternary
$(\sigma,\tau,\xi)-$derivation such that $\sigma,\tau$ and $\xi$ are
linear mappings on $A$, additionally, $\tau(e)=e$, then it is easy
to prove that $D$ is a Lie $(\sigma,\xi)-$derivation.
\\

The stability of functional equations was started in 1940 with a
problem raised by S. M. Ulam [19].
 In 1941 Hyers affirmatively solved the problem of S. M. Ulam in the context of Banach spaces.
 In 1950 T.Aoki [4] extended the Hyers' theorem. In 1978, Th. M. Rassias [16] formulated and proved the
 following Theorem:\\
Assume that $E_1$ and $E_2$ are real normed spaces with $E_2$ complete, $f:E_1\to E_2$ is a mapping such that
for each fixed $x\in E_1$ the mapping $t\to f(tx)$ is continuous on $\Bbb R$, and let there exist $\epsilon\geq 0$
 and $p\in [0,1)$ such that $\|f(x+y)-f(x)-f(y)\|\leq\epsilon(\|x\|^p+\|y\|^p)$ for all $x,y\in E_1$.
 Then there exists a unique linear mapping $T:E_1\in E_2$ such that $\|f(x)-T(x)\|\leq\epsilon\frak{\|x\|^p}{(1-2^p)}$
 for all $x\in E_1.$\\
The equality $\|f(x+y)-f(x)-f(y)\|\leq\epsilon(\|x\|^p+\|y\|^p$ has
provided extensive influence in the development of of what we now
call Hyers--Ulam--Rassias stability of functional equations
[8,11,15,17,18]. In 1994, a generalization of Rassias' theorem was
obtained by Gavruta [9], in which he replaced the bound
$\epsilon(\|x\|^p+\|y\|^p)$ by a general control function. \vskip
5mm
\section{\bf Lie ternary $({\sigma},{\tau},{\xi})-$derivations on Banach ternary algebras}
\vskip 2mm

In this section our aim is to establish the Hyers--Ulam--Rassias stability of Lie ternary
$(\sigma,\tau,\xi)-$derivations.\\
\begin{thm}
Suppose $f:A \to X$ is a mapping with $f(0)=0$ for which there exist mappings $g, h, k:A\to A$ with $g(0)=h(0)=k(0)=0$ and a function $\varphi:A\times A\times A\times A\times A\to [0,\infty]$ such that
$$\widetilde{\varphi}(x,y,u,v,w)=\frac{1}{2}\sum_{n=0}^{\infty}\varphi(2^n x,2^n y,2^n u,2^n v,2^n w)<\infty\eqno(2.2)$$
\begin{align*}
&\|f(\lambda x+\lambda y+[uvw]_A)-\lambda f(x)-\lambda f(y)-[[f(u)vw]_X]_{(g,h,k)}+[[f(v)uw]_X]_{(g,h,k)}+\\
&[[f(w)vu]_X]_{(g,h,k)}\|\leq \varphi(x,y,u,v,w)\hspace{7.2 cm}(2.3)
\end{align*}
$$\|g(\lambda x+\lambda y)-\lambda g(x)-\lambda g(y)\|\leq\varphi(x,y,0,0,0)$$
$$\|h(\lambda x+\lambda y)-\lambda h(x)-\lambda h(y)\|\leq\varphi(x,y,0,0,0)$$
$$\|k(\lambda x+\lambda y)-\lambda k(x)-\lambda k(y)\|\leq\varphi(x,y,0,0,0)$$
\\
for all $\lambda\in \Bbb T^1(:=\{\lambda \in \Bbb C~;
|\lambda|=1\})$ and for all $x,y,u,v,w\in A.$ Then there exist
unique linear mappings $\sigma, \tau$ and $\xi$ from A to A
satisfying $$\|g(x)-\sigma(x)\|\leq
\widetilde{\varphi}(x,x,0,0,0)\eqno(2.4)$$
$$\|h(x)-\tau(x)\|\leq \widetilde{\varphi}(x,x,0,0,0)\eqno(2.5)$$
and$$\|k(x)-\xi(x)\|\leq
\widetilde{\varphi}(x,x,0,0,0)\eqno(2.6)$$and there exist a unique
Lie ternary$(\sigma,\tau,\xi)-$derivation on $D:A\to X$ such that
 $$\|f(x)-D(x)\|\leq\widetilde{\varphi}(x,x,0,0,0)\eqno(2.7)$$
 for all $x\in A.$
\end{thm}
\begin{proof}
One can show that the limits
$$\sigma(x):=\lim_n\frac{1}{2^n}{g(2^nx)}$$
$$\tau(x):=\lim_n\frac{1}{2^n}{h(2^nx)}$$
$$\xi(x):=\lim_n\frac{1}{2^n}{k(2^nx)}$$
exist for all $x\in A$, also $\sigma,\tau$ and  $\xi$  are unique
linear mappings which satisfy (2.4), (2.5) and  (2.6) respectively
(see [17]).

Put $\lambda=1$ and $u=v=w=0$ in (2.3) to obtain
$$\|f(x+y)-f(x)-f(y)\|\leq \phi(x,y,0,0,0)\hspace{1.5
cm}(x,y\in A).\eqno(2.8)$$ Fix $x\in A.$ Replace $y$ by $x$ in (2.8)
to get
$$\|f(2x)-2f(x)\|\leq\varphi(x,x,0,0,0).$$
One can use the induction to show that
$$\|\frac{f(2^px)}{2^p}-\frac{f(2^qx)}{2^q}\|\leq\frac{1}{2}\sum_{k=q}^{p-1}\varphi(2^kx,2^kx,0,0,0)\eqno(2.9)$$
for all $x\in A,$ and all $p>q\geq0.$ It follows from the
convergence of series (2.2) that the sequence
$\{\frac{f(2^nx)}{2^n}\}$ is Cauchy. By the completeness of $X$,
this sequence is convergent. Set
$$D(x)=\lim_{n\to \infty}\frac{f(2^nx)}{2^n}\eqno(2.10)$$
for all $x\in A$. Putting $u=v=w=0$ and replacing $x,y$ by $2^nx$
and $2^ny$ in (2.3) respectively,
 and divide the both sides of the inequality by $2^n$ we get
$$\|2^{-n}f(2^n(\lambda x+\lambda y))-2^{-n}\lambda f(2^nx)-2^{-n}\lambda f(2^ny)\|\leq\frac{1}{2^n}\varphi(2^nx,2^nx,0,0,0).$$
Passing to the limit as $n\to\infty$ we obtain $D(\lambda x+\lambda y)=\lambda D(x)+\lambda D(y).$\\
Put $q=0$ in (2.9) to get
$$\|\frac{f(2^px)}{2^p}-f(x)\|\leq\frac{1}{2}\sum_{k=0}^{p-1}\varphi(2^kx,2^kx,0,0,0)$$
for all $x\in A.$\\
Taking the limit as $p\to\infty$ we infer that
$$\|f(x)-D(x)\|\leq\widetilde{\varphi}(x,x,0,0,0)$$
for all $x\in A.$ Next, let $\gamma\in\Bbb C (\gamma\neq0)$ and let
$N$ be a positive integer number greater than $|\gamma|.$ It is
shown that there exist two numbers $\lambda_1,\lambda_2\in \Bbb T$
such that $2\frac{\gamma}{N}=\lambda_1+\lambda_2.$ since $D$ is a
additive, we have $D(\frac{1}{2}x)=\frac{1}{2}D(x)$ for all $x\in
A.$ Hence
\begin{align*}
D(\gamma x)&=D(\frac{N}{2}.2.\frac{\gamma}{N}x)=ND(\frac{1}{2}.2.\frac{\gamma}{N}x)=\frac{N}{2}D(2.\frac{\gamma}{N}x)\\
&=\frac{N}{2}D(\lambda_1x+\lambda_2x)=\frac{N}{2}(D(\lambda_1x)+D(\lambda_2x))\\
&=\frac{N}{2}(\lambda_1+\lambda_2)D(x)=(\frac{N}{2}.2.\frac{\gamma}{N})D(x)=\gamma D(x)
\end{align*}
for all $x\in A.$ Thus $D$ is linear.\\
Suppose that there exists another ternary
$(\sigma,\tau,\xi)-$derivation $D^{'}:A \to X$ satisfying $(2.7).$
Since $D^{'}(x)=\frac{1}{2^n}D^{'}(2^n x),$ we see that
\begin{align*}
&\|D(x)-D^{'}(x)\|=\frac{1}{2^n}\|D(2^n
x)-D^{'}(2^n x)\|\\
&\leq \frac{1}{2^n}(\|f(2^n x)-D(2^n x)\|+\|f(2^n x)-D^{'}(2^n x)\|)\\
&\leq 4\theta \frac{2^p}{2-2^p}2^{n(p-1)}\|x\|^p~,
\end{align*}
which tends to zero as $n \to \infty$ for all $x \in A.$ Therefore
$D^{'}=D$ as claimed. Similarly one can use (2.4), (2.5) and (2.6) to show that there exist unique
linear mappings $\sigma, \tau$ and $\xi$ defined
by $\lim_{n\to \infty} \frac{g(2^nx)}{2^n},\lim_{n\to \infty}\frac{h(2^nx)}{2^n}$
and $\lim_{n\to \infty}\frac{k(2^nx)}{2^n}$, respectively.\\
 Putting $x=y=0$ and replacing $u,v,w$ by $2^nu, 2^nv$ and $2^nw$ in (2.3) respectively, we obtain
 \begin{align*}
 &\|f([2^{3n}uvw]_A)-[[f(2^nu)2^{2n}vw]_X]_{(g,h,k)}+[[f(2^nv)2^{2n}uw]_X]_{(g,h,k)}+[[f(2^nw)2^{2n}vu]_X]_{(g,h,k)}\|\\
 &\leq\varphi(0,0,2^nu,2^nv,2^nw),
 \end{align*}
 then
 \begin{align*}
 &\frac{1}{2^{3n}}\|f([2^{3n}uvw]_A)-[[f(2^nu)2^{2n}vw]_X]_{(g,h,k)}+[[f(2^nv)2^{2n}uw]_X]_{(g,h,k)}+[[f(2^nw)2^{2n}vu]_X]_{(g,h,k)}\|\\
 &\leq\frac{1}{2^{3n}}\varphi(0,0,2^nu,2^nv,2^nw)
 \end{align*}
 for all $u,v,w\in A,$ hence,
 \begin{align*}
 &\lim_{n\to\infty}\frac{1}{2^{3n}}\|f([2^{3n}uvw]_A)-[[f(2^nu)2^{2n}vw]_X]_{(g,h,k)}+[[f(2^nv)2^{2n}uw]_X]_{(g,h,k)}\\
 &+[[f(2^nw)2^{2n}vu]_X]_{(g,h,k)}\|\leq\lim_{n\to\infty}\frac{1}{2^{3n}}\varphi(0,0,2^nu,2^nv,2^nw)\\
 &=0
 \end{align*}
 therefore
 \begin{align*}
D&([uvw]_A)=\lim_{n\to\infty}\frac{f(2^{3n}[uvw]_A)}{2^{3n}}=\lim_{n\to\infty}\frac{f([2^nu2^nv2^nw]_A)}{2^{3n}}\\
 &=\lim_{n\to\infty}(\frac{[[f(2^nu)2^nv2^nw]_X]_{(g,h,k)}-[[f(2^nv)2^nu2^nw]_X]_{(g,h,k)}-[[f(2^nw)2^nv2^nu]_X]_{(g,h,k)}}{2^{3n}})\\
 &=\lim_{n\to\infty}(\frac{f(2^nu)h(2^nv)k(2^nw)-g(2^nw)h(2^nv)f(2^nu)-f(2^nv)h(2^nu)k(2^nw)}{2^{3n}}\\
&+\frac{g(2^nw)h(2^nu)f(2^nv)-f(2^nw)h(2^nv)k(2^nu)+g(2^nu)h(2^nv)f(2^nw)}{2^{3n}})\\
&=(D(u)\tau(v)\xi(w)-\sigma(w)\tau(v)D(u))-(D(v)\tau(u)\xi(w)-\sigma(w)\tau(u)D(v))\\
& -(D(w)\tau(v)\xi(u)-\sigma(u)\tau(v)D(w))\\
&=[[D(u)vw]_X]_{(\sigma,\tau,\xi)}-[[D(v)uw]_X]_{(\sigma,\tau,\xi)}-[[D(w)vu]_X]_{(\sigma,\tau,\xi)}
 \end{align*}
 for each $u,v,w\in A.$ Hence, the linear mapping $D$ is a Lie ternary $(\sigma,\tau,\xi)-$derivation.
\end{proof}
\begin{cor}Suppose $f:A\to X$ is a mapping with $f(0)=0$ for which there exist mappings $g, h, k:A\to A$
with $g(0)=h(0)=k(0)=0$ and there exists
$\theta\geq0$ and $p\in[0,1)$ such that

 \begin{align*}
&\|f(\lambda x+\lambda y+[uvw]_A)-\lambda f(x)-\lambda f(y)-[[f(u)vw]_X]_{(g,h,k)}+[[f(v)uw]_X]_{(g,h,k)}\\
&+[[f(w)vu]_X]_{(g,h,k)}\|\leq\theta(\|x\|^p+\|y\|^p+\|u\|^p+\|v\|^p+\|w\|^p),
 \end{align*}
 \\
$$\|g(\lambda x+\lambda y)-\lambda g(x)-\lambda g(y)\|\leq\theta(\|x\|^p+\|y\|^p)$$
$$\|h(\lambda x+\lambda y)-\lambda h(x)-\lambda h(y)\|\leq\theta(\|x\|^p+\|y\|^p)$$
$$\|k(\lambda x+\lambda y)-\lambda k(x)-\lambda k(y)\|\leq\theta(\|x\|^p+\|y\|^p)$$
\\
for all $\lambda\in \Bbb T=\{\lambda\in \Bbb C : |\lambda|=1\}$ and
for all $x,y\in A.$ Then there exist unique linear mappings $\sigma,
\tau$ and $\xi$ from $A$ to $A$ satisfying
$\|g(x)-\sigma(x)\|\leq\frac{\theta\|x\|^p}{1-2^{p-1}},
\|h(x)-\tau(x)\|\leq\frac{\theta\|x\|^p}{1-2^{p-1}}$ and
$\|k(x)-\xi(x)\|\leq\frac{\theta\|x\|^p}{1-2^{p-1}}$, and there
exists a unique Lie ternary $(\sigma, \tau, \xi)-$derivation $D:A\to
X$ such that
$$\|f(x)-D(x)\|\leq\frac{\theta\|x\|^p}{1-2^{p-1}}\eqno(2.11)$$
for all $x\in A.$
\end{cor}
\begin{proof}
Put $\varphi(x,y,u,v,w)=\theta(\|x\|^p+\|y\|^p+\|u\|^p+\|v\|^p+\|w\|^p)$ in Theorem 2.1.
\end{proof}

Note that a linear mapping $D:(A,[~]_A) \to (X,[~]_X)$ is called a
Jordan Lie ternary $(\sigma,\tau,\xi)-$derivation, if
$$D([aaa]_A)=[[D(a)aa]_X]_{(\sigma,\tau,\xi)}+[[D(a)aa]_X]_{(\sigma,\tau,\xi)}+[[D(a)aa]_X]_{(\sigma,\tau,\xi)}$$
 for all
$a\in A$,

\begin{thm}
Suppose $f:A \to X$ is a mapping with $f(0)=0$ for which there exist
mappings $g, h, k:A\to A$ with $g(0)=h(0)=k(0)=0$ and a function
$\varphi:A\times A\times A\to [0,\infty]$ such that
$$\widetilde{\varphi}(x,y,u)=\frac{1}{2}\sum_{n=0}^{\infty}\varphi(2^n x,2^n y,2^n u)<\infty\eqno(2.12)$$
\begin{align*}
&\|f(\lambda x+\lambda y+[uuu]_A)-\lambda f(x)-\lambda f(y)-[[f(u)uu]_X]_{(g,h,k)}+[[f(u)uu]_X]_{(g,h,k)}+\\
&[[f(u)uu]_X]_{(g,h,k)}\|\leq \varphi(x,y,u)\hspace{7.8 cm}(2.13)
\end{align*}
$$\|g(\lambda x+\lambda y)-\lambda g(x)-\lambda g(y)\|\leq\varphi(x,y,0)$$
$$\|h(\lambda x+\lambda y)-\lambda h(x)-\lambda h(y)\|\leq\varphi(x,y,0)$$
$$\|k(\lambda x+\lambda y)-\lambda k(x)-\lambda k(y)\|\leq\varphi(x,y,0)$$
\\
for all $\lambda\in \Bbb T^1(:=\{\lambda \in \Bbb C~;
|\lambda|=1\})$ and for all $x,y,u\in A.$ Then there exist unique
linear mappings $\sigma, \tau$ and $\xi$ from A to A, and a unique
Jordan Lie ternary $(\sigma,\tau,\xi)-$derivation $D:A \to X$
satisfying (2.4), (2.5), (2.6) and (2.7), respectively.
\end{thm}
\begin{proof}
By the same reasoning as the proof of Theorem 2.1, the limits
$$D(x):=\lim_n\frac{1}{2^n}{f(2^nx)}$$
$$\sigma(x):=\lim_n\frac{1}{2^n}{g(2^nx)}$$
$$\tau(x):=\lim_n\frac{1}{2^n}{h(2^nx)}$$
$$\xi(x):=\lim_n\frac{1}{2^n}{k(2^nx)}$$
exist for all $x\in A$, also $\sigma,\tau,\xi$ and $D$  are unique
linear mappings which satisfy (2.4), (2.5), (2.6) and (2.7)
respectively.
 Putting $x=y=0$ and replacing $u$ by $2^nu$  in (2.13), we obtain

 \begin{align*}
\|D([uuu]_A)&-[[D(u)uu]_X]_{(\sigma,\tau,\xi)}-[[D(u)uu]_X]_{(\sigma,\tau,\xi)}-[[D(u)uu]_X]_{(\sigma,\tau,\xi)}\|\\
&=\lim_{n\to\infty}\|\frac{f(2^{3n}[uuu]_A)}{2^{3n}}-(D(u)\tau(u)\xi(u)-\sigma(u)\tau(u)D(u))\\
&-(D(u)\tau(u)\xi(u) -\sigma(u)\tau(u)D(u)) -(D(u)\tau(u)\xi(u)-\sigma(u)\tau(u)D(u))\|\\
&=\lim_{n\to\infty}\|\frac{f([2^nu2^nu2^nu]_A)}{2^{3n}}\\
&-
(\frac{f(2^nu)h(2^nu)k(2^nu)-g(2^nu)h(2^nu)f(2^nu)-f(2^nu)h(2^nu)k(2^nu)}{2^{3n}}\\
&+\frac{g(2^nu)h(2^nu)f(2^nu)-f(2^nu)h(2^nu)k(2^nu)+g(2^nu)h(2^nu)f(2^nu)}{2^{3n}})\|\\
&=\lim_{n\to\infty}\|\frac{f([2^nu2^nu2^nu]_A)}{2^{3n}}\\
&-
(\frac{[[f(2^nu)2^nu2^nu]_X]_{(g,h,k)}-[[f(2^nu)2^nu2^nu]_X]_{(g,h,k)}-[[f(2^nu)2^nu2^nu]_X]_
 {(g,h,k)}}{2^{3n}})\|\\
 &\leq\lim_{n\to\infty}\frac{1}{2^{3n}}\varphi(0,0,2^nu)\\
 &=0
 \end{align*}
 for each $u\in A.$ Hence, the linear mapping $D$ is a Jordan Lie ternary $(\sigma,\tau,\xi)-$derivation.
\end{proof}
\begin{cor}Suppose $f:A\to X$ is a mapping with $f(0)=0$ for which there exist mappings $g, h, k:A\to A$
with $g(0)=h(0)=k(0)=0$ and there exists $\theta\geq0$ and
$p\in[0,1)$ such that

 \begin{align*}
&\|f(\lambda x+\lambda y+[uuu]_A)-\lambda f(x)-\lambda f(y)-[[f(u)uu]_X]_{(g,h,k)}+[[f(u)uu]_X]_{(g,h,k)}\\
&+[[f(u)uu]_X]_{(g,h,k)}\|\leq\theta(\|x\|^p+\|y\|^p+\|u\|^p),
 \end{align*}
 \\
$$\|g(\lambda x+\lambda y)-\lambda g(x)-\lambda g(y)\|\leq\theta(\|x\|^p+\|y\|^p)$$
$$\|h(\lambda x+\lambda y)-\lambda h(x)-\lambda h(y)\|\leq\theta(\|x\|^p+\|y\|^p)$$
$$\|k(\lambda x+\lambda y)-\lambda k(x)-\lambda k(y)\|\leq\theta(\|x\|^p+\|y\|^p)$$
\\
for all $\lambda\in \Bbb T=\{\lambda\in \Bbb C : |\lambda|=1\}$ and
for all $x,y\in A.$ Then there exist unique linear mappings $\sigma,
\tau$ and $\xi$ from $A$ to $A$ satisfying
$\|g(x)-\sigma(x)\|\leq\frac{\theta\|x\|^p}{1-2^{p-1}},
\|h(x)-\tau(x)\|\leq\frac{\theta\|x\|^p}{1-2^{p-1}}$ and
$\|k(x)-\xi(x)\|\leq\frac{\theta\|x\|^p}{1-2^{p-1}}$, and there
exists a unique Jordan Lie ternary $(\sigma, \tau, \xi)-$derivation
$D:A\to X$ such that
$$\|f(x)-D(x)\|\leq\frac{\theta\|x\|^p}{1-2^{p-1}}\eqno(2.11)$$
for all $x\in A.$
\end{cor}
\begin{proof}
Put $\varphi(x,y,u)=\theta(\|x\|^p+\|y\|^p+\|u\|^p)$ in Theorem 2.3.
\end{proof}

{\small

}
\end{document}